\documentstyle[11pt]{article}
\newtheorem {th}{Theorem}
\newtheorem {lem}[th]{Lemma}

\newtheorem {cor}[th]{Corollary}

\makeatletter
\renewcommand \thesection {\@arabic\c@section.}
\makeatother
\def\sect#1{\section[#1]{\hskip -10pt
#1.}}

\def\eqalign#1{\begin{array}{rl} #1 \end{array}}

\def\Cox{\hfill \Box}
\def\DD{\Delta}
\def\DE{\Delta'}

\def\E{{\bf{E}}}

\def\P{{\bf{P}}}

\def\F{{\cal{F}}}
\def\L{{\cal{L}}}
\def\|{\, | \, }
\def\one{{\bf 1}}
\def\itemm{\par\hang\textindent}
\def\IN{{\bf N}}
\def\real{{\Re e }}

\begin{document}

\begin{titlepage}
\begin{center}
{\large \bf Moment conditions for a sequence with negative drift \\
to be uniformly bounded in $L^r$} \\
\end{center}
\vspace{1ex}
\begin{center}
Robin Pemantle\footnote{Department 
of Mathematics, University of Wisconsin--Madison, Van Vleck Hall, 480 Lincoln
Drive, Madison, WI 53706.  Internet: {\tt pemantle@math.wisc.edu}.
Research supported in part by 
NSF grants DMS-9300191 and DMS-9803249.}
\quad and \quad
Je{f}frey S.\ Rosenthal\footnote{Department of Statistics, University of
Toronto, Toronto, Ontario, Canada M5S 3G3.  
Internet: {\tt jeff@utstat.toronto.edu}.
Research supported in part by NSERC of Canada.}
\end{center}

\vspace{2ex}
\centerline{(August 27, 1998; revised January 12, 1999.)}

\vfill

\noindent {\bf ABSTRACT:} \hfill \break
Suppose a sequence of random variables $\{ X_n \}$ has negative drift
when above a certain threshold and has increments bounded in $L^p$.
When $p > 2$ this implies that $\E X_n$ is bounded above by a constant
independent of $n$ and the particular sequence $\{ X_n \}$.  When $p \leq
2$ there are counterexamples showing this does not hold.  In general,
increments bounded in $L^p$ lead to a uniform $L^r$ bound on $X_n^+$
for any $r < p-1$, but not for $r \ge p-1$.  These results are motivated
by questions about stability of queueing networks.

\vfill



\bigskip

\end{titlepage}

\setcounter{equation}{0}

\sect{Introduction}

Let $X_0,X_1,\ldots$ be a sequence of real-valued random variables. 
We wish to find a condition, along the lines of behaving like a 
supermartingale when sufficiently large, that will guarantee
$\sup_n \E X_n < \infty$.  In particular, we do not wish to assume
any special properties of the increments such as independence, 
$r$-dependence, Markov property, symmetry, discreteness or nondiscreteness.
Under what conditions can we guarantee that $\sup_n \E(X_n) < \infty$?  

Specifically, we suppose that for some $a>0$ and some $J$ we have
$$
\E (X_{n+1} - X_n \| X_0 , \ldots , X_n) \leq -a \quad \mbox{ on the event }
\quad \{X_n > J\}
\eqno(\mbox{C1})
$$
for all $n$.  That is, the process has {\it negative drift} whenever it
is above the point $J$.  This condition alone says nothing about 
possible large jumps out of the interval $(-\infty , J)$, so we also
assume that for some $p \ge 1$ and some $V<\infty$ we have
$$
\E (|X_{n+1} - X_n|^p \| X_0 , \ldots , X_n) \leq V
\eqno(\mbox{C2})
$$
for all $n$.  That is, the process has {\it increments
with bounded $p^{\rm th}$ moments}.  Conditions (C1) and (C2) are
meant to characterize the behavior of sequences attracted to a basin,
which always decreases in expectation except when it is already small.
(One such example is a nonnegative Lyapunov function of a Markov chain;
see for example Meyn and Tweedie, 1993.)
Do these two conditions together imply that $\sup_n \E(X_n) < \infty$?

When $p=1$ the answer is no.  Although an honest supermartingale
has $\E X_n \leq \E X_0$ for all $n$, a process behaving like a 
supermartingale above $J$ may have $\E X_n$ unbounded.  More 
surprisingly, the answer is still no for $p = 2$.  However, if
one assumes (C1) and (C2) with $p > 2$, then necessarily
$\sup_n \E X_n < \infty$.  More generally, we show
(Theorem~\ref{th:1}) that (C1) and (C2) imply that $\sup_n \E\left(
(X_n^+)^r \right) < \infty$ whenever $r < p-1$, and that this
bound on $r$ in terms of $p$ (or $p$ in terms of $r$) is sharp.  
Furthermore, our bounds are may be explicitly computed and depend 
only on the parameters $a$, $J$, $p$, $V$, and $r$.

Our results are motivated by questions about queueing networks.
Specifically, several authors (Borodin et al., 1996, 1998; Andrews
et al., 1996) consider network loads under the influence of a {\it
stochastic adversary}.  Here $X_n$ is the load of the network at time $n$.
The adversary may add new packets to the network in virtually any manner,
subject only to a load condition which leads to (C1) plus a moment condition 
such as (C2).  (The load condition corresponds to the statement that, once 
the network is operating at full capacity, it processes packets more quickly
on average than the adversary can add them.)  The network is considered
to be {\it stable} if the expected load remains bounded, i.e.\ if $\sup_n
\E(X_n) < \infty$.  In this context, our Corollary~\ref{cor:X_0} may
be interpreted as saying that a queueing network in the presence of a
stochastic adversary is guaranteed to be stable, provided it satisfies
the load condition (C1), and also the moment condition (C2) for some $p>2$.
On the other hand, if $p \le 2$ then there is no such guarantee.

We note that there has been some previous work on related questions.
For example, Hastad et al.\ (1996) consider bounds on $\sup_n \E(X_n)$
for certain time-homogeneous {\it Markovian} systems which correspond
to particular ``backoff protocols'' for resolving ethernet conflicts.
Close to our work, Hajek (1982) investigates bounds on hitting times
for general random sequences having bounded {\it exponential} moments,
and derives corresponding bounds on exponential moments of the hitting
times; his work may thus be seen, roughly, as the $p\to\infty$ limit of
our bounds.

Finally, we note that while the notion of ``stability'' considered here
(namely, that $\sup_n \E(X_n) < \infty$) is different from that of
Markov chain stability (see e.g.\ Meyn and Tweedie, 1993), there are
some connections.  For example, it is known (see Tweedie, 1983, Theorem
2) that for $k\in\IN$,
if $\{X_n\}$ is an aperiodic Harris-recurrent time-homogeneous
Markov chain having stationary distribution $\pi(\cdot)$, and if $m_k \equiv
\int x^k \, \pi(dx) < \infty$, then for $\pi$-a.e.\ $x$, $\E_{\delta_x}
(X_n^k) \to m$, and hence $\sup_n \E_{\delta_x} (X_n^k) < \infty$.  In other
words, for such a Markov chain, stability in our sense is implied by
standard Markov chain stability.  In fact, it is known (e.g.\ Tuominen
and Tweedie, 1994) that when $\{ X_n \}$ is a random walk with negative
drift, reset to zero when it attempts to leave the nonnegative half-line
and having square integrable increments, then $\E_{\delta_x} (X_n)$
will converge and hence be bounded.  This shows that our (C1) and (C2)
do represent a greater generality than the random walk context.

The paper is organized as follows.  In Section~2 we state
the main result, along with two extensions.  (The extensions are
reasonably straightforward, but we include them in order to provide
readily referenceable results that don't assume more than is needed.)
We also provide in Section~2 a simpler proof of the main theorem
in the case where $p > 4$ and $r = 1$,
since in this case the back-of-the-napkin computation works, and
anyone not interested in the sharp moment condition need read no further.
In Section~3, we give examples to show why (C2) is needed with $p - 1 > r$
and why it is important to have moment bounds for the negative part 
of the increments as well as the positive part.  Proofs are given
in Sections~4 and~5, with Section~4 containing a reduction to a
result on martingales and Section~5 containing a proof of the
martingale result.

\sect{Main results}

Throughout this paper, the filtration $\{ \F_n \}$ refers to any 
filtration to which $\{ X_n \}$ is adapted.  We continue to use (C1)
and (C2) for conditional expectations with respect to $\F_n$, slightly
generalizing the notation of the introductory section.  

\begin{th} \label{th:1}
Let $X_n$ be random variables and suppose that there exist constants
$a>0$, $J$, $V<\infty$, and $p > 2$, such that $X_0 \le J$, and for all $n$,
$$
\E (X_{n+1} - X_n \| \F_n) \leq -a \quad \mbox{ on the event }
\quad \{X_n > J\}
\eqno(\mbox{\rm C1})
$$
and
$$
\E (|X_{n+1} - X_n|^p \| X_0 , \ldots , X_n) \leq V
\eqno(\mbox{\rm C2})
$$
Then for any $r \in (0 , p-1)$ there is a $c = c (p , a, V, J , r) > 0$ 
such that $\E (X_n^+)^r < c$ for all $n$.
\end{th}

Applying this theorem to the process $X_n' := X_n - (X_0-J)^+$ in the
case $r = 1$ immediately yields the corollary:

\begin{cor} \label{cor:X_0}
Under hypotheses (C1) and (C2) of Theorem~\ref{th:1}, 
but without assuming $X_0 \le J$, we have
$$\E (X_n \| \F_0) \leq c(p,a,V,J,1) + (X_0-J)^+ .$$
\end{cor}

\bigskip
\noindent \bf Remark. \rm
By following through the proof (presented in Sections 4 and 5), we are able
to provide an explicit formula for the quantity $c(p,a,V,J,r)$ of Theorem 1.  
Indeed, for $a=1$ and $J=0$, we have
$c(p,1,V,0,r) = K \zeta(p-r)$
where
$K = C(b,p,r) = 2^{p/2} c_2 C'(p,b) + c_4$.  Here
$b = 2^p (B + (1+B)^p)$;
$B = 2^p (1+V)$;
$C'(p,b) = \max(1, \ c'(p,b))$;
$c'(p,b) = c_p b (1+c_p^{-1})^p$;
$c_2 = c_p b (4^p + 4^{p-r} {r \over p-r})$;
$c_4 = C'(p,b) c_3 \zeta(p/2)$;
$c_3 = 3^r 4^p b (c_p b + {p \over p-r} + 3^r)$;
and $c_p = (p-1)^p$ is the constant from
Burkholder's inequality.
(Recall that $\zeta(w) \equiv \sum_{i=1}^\infty i^{-w}$ is the
Riemann zeta function, finite for $\real(w) > 1$.)
Then for general $a$ and $J$, we have
$c(p,a,V,J,r) = J + a^r c(p,1,V/a^p,0,r)$.
Now, these formulae are clearly rather messy, and may be of limited
practical use.  However, it may still be helpful to have them
available for ready reference.
\bigskip

We also state an extension allowing the negative part of the increments
to avoid the moment condition in (C2):

\begin{cor} \label{th:9}
The conclusion of Theorem~\ref{th:1} still holds when 
$X_{n+1} - X_n$ is replaced by $(X_{n+1} - X_n) \one_{X_{n+1} - X_n > Z_n}$
in conditions (C1) and (C2), and $Z_n \leq -a$ is any sequence adapted to 
$\{ \F_n \}$.
\end{cor}

The proof of Theorem~\ref{th:1} proceeds by decomposing according
to the last time $U$ before time $n$ that $\{ X_k \}$ was less than $J$.
When $p > 4$, Markov's inequality, together with a crude $L^p$ estimate
on $X_n - X_U$, gives bounds on the tails of $X_n$ sufficient
to yield Theorem~\ref{th:1}.  We finish the section by giving this argument.

Assume the notation and hypotheses of Theorem 1.  Fix a positive
integer $n$.  Let $U = \max\{ k \le n; \ X_k \le J \}$.
Let $\mu_i = \E\left( X_{i+1}-X_i \, | \, \F_i \right)$, so
that $\mu_i \le -a$ on $\{X_i>J\}$.  We may recenter (see the
proof of Corollary~\ref{cor:1} for details) to obtain
$$
\E\left( (X_{n+1}-X_n-\mu_n)^p \, | \, \F_n \right)
\ \le \ V'
$$
for some $V' < \infty$.

But then, for $t > J$, we have
$$
\P( X_n \ge t)
= \sum_{k=0}^{n-1} \P(X_n \ge t, \ U=k)
\qquad\qquad\qquad
$$
$$
\le \sum_{k=0}^{n-1} \P(X_n - X_k \ge t-J, \ X_k \le J, \ X_i > J \ {\rm for} \
k<i<n)
$$
$$
\le \sum_{k=0}^{n-1} \P\bigg( (X_n - X_{n-1}-\mu_{n-1}) + \ldots 
+ (X_{k+1}-X_k-\mu_{k}) \ge t-J-V^{1/p}+a(n-k-1), 
$$ 
$$ X_k \le J, \ X_i > J \ {\rm for} 
\ k<i<n \bigg)
$$
[since $\mu_i \le -a$ for $k<i<n$, and $\mu_k \le V^{1/p}$]
$$
\le \sum_{k=0}^{n-1} \E\bigg( \Big| (X_n - X_{n-1}-\mu_{n-1}) + \ldots 
+ (X_{k+1}-X_k-\mu_{k}) \Big|^p \bigg) \Big( t-J-V^{1/p}+a(n-k-1) \Big)^{-p}
$$
[by Markov's inequality]
$$
\le \sum_{k=0}^{n-1} c_p V' (n-k)^{p/2} \Big( t-J-V^{1/p}+a(n-k-1) \Big)^{-p}
$$
[by Lemma \ref{lem:Lp}, which is a direct application of Burkholder's
inequality]
$$
\le \sum_{\ell=0}^{\infty} c_p V' (\ell+1)^{p/2}
\Big( t-J-V^{1/p}+a\ell \Big)^{-p} \, .
$$

It then follows that
$$
\E(X_n) \ = \
\int_0^\infty dt \, \P(X_n \ge t)
\qquad\qquad\qquad
$$
$$
\le \
(J+V^{1/p}+1) \ + \ \int_{J+V^{1/p}+1}^\infty dt \, 
\sum_{\ell=0}^{\infty} c_p V' (\ell+1)^{p/2}
\left( t-J-V^{1/p}+a\ell \right)^{-p}
$$
This integral-of-sum does not depend on $n$.  Furthermore, for $p>4$ it is
straightforward to check (by integrating first) that it
is finite.  This gives the result.  $\Cox$

\sect{Some counterexamples}

We here present a few counterexamples to show that the hypotheses of
Theorem 1 (in particular, the restriction that $p>2$) are really
necessary.

\itemm{1.}
The following example is due to Madhu Sudan (personal communication
via A.~Borodin).
Let $\{X_n\}$ be a time-inhomogeneous Markov chain 
such that $\P(X_2=0) = 2/3 = 1 - \P(X_2=2)$, and such that
for $n \ge 2$,
$$
        \P(X_{n+1} = n+1 \, | \, X_n = n) = 1 - 2/n
$$
$$
        \P(X_{n+1} = 0 \, | \, X_n = n) = 2/n
$$
$$
        \P(X_{n+1} = n+1 \, | \, X_n = 0) = 1/n
$$
$$
        \P(X_{n+1} = 0 \, | \, X_n = 0) = 1 - 1/n
$$
These transition probabilities were chosen to ensure that
$$
X_n \ = \
\left\{
\eqalign{
0,& \quad {\rm prob} \ 2/3 \cr
n,& \quad {\rm prob} \ 1/3 }
\right.
$$
for all $n \ge 2$.  Hence, $\E(X_n) = n/3$, so that
$\sup_n \E(X_n) = \infty$.

\itemm{}
On the other hand, it is easily verified
that (C1) is satisfied with $a=1$ and $J=0$.
Furthermore, (C2) is satisfied with $p=1$ and $V=3$.
We conclude that condition (C1) alone, or combined with
(C2) with $p=1$, does not guarantee stability.

\itemm{2.}
When (C2) holds with $p=2$ it appears one has to do a little
more to engineer a counterexample; specifically, we line up all 
the jumps out of $(-\infty , J)$ to amass at a fixed time $M$.
Fix a large integer $M$, and define a time-inhomogeneous Markov chain 
by setting $X_0=0$, and, for $0 \le n \le M-1$, letting
$$
        \P(X_{n+1} = X_n-1 \, | \, X_n > 0) = 1
$$
$$
        \P(X_{n+1} = 0 \, | \, X_n = 0 ) = 1 - (M-n)^{-2}
$$
$$
        \P(X_{n+1} = 2(M-n) \, | \, X_n = 0 ) = (M-n)^{-2}
$$
Then it is easily verified that
(C1) is again satisfied with $a=1$ and $J=0$.
Furthermore, (C2) is satisfied with $p=2$ and $V=4$.

\itemm{}
On the other hand, 
setting $A = \exp\left( - \sum_{i=1}^\infty 1/i^2 \right) > 0$,
we compute that
$$ \E(X_M) =
\sum_{k=0}^M (M-k+1) \, \P (X_k > 0 \ {\rm and} \
 X_j = 0 \ {\rm for} \ j < k)
$$
$$
\ge \ \sum_{k=0}^{M} (M-k+1) \left( A / (M-k+1)^2 \right)
$$
$$
                = \ \sum_{k=0}^{M} A / (M-k+1)
                \ = \ \sum_{j=1}^{M+1} A / j
$$
which goes to infinity (like $A \log M$) as $M \to \infty$.

\itemm{}
This shows that $\E X_n$ cannot be bounded in terms of $a, J$ and $V$, and
by ``stringing together'' such examples, for larger and
larger choices of $M$, we can clearly make $\sup_n \E(X_n) = \infty$.
We conclude that condition (C1), combined with
(C2) with $p=2$, still does not guarantee stability of $\{X_n\}$.

\itemm{3.}
>From the queueing theory perspective, it would be desirable,
in condition (C2) of Theorem~\ref{th:1}, to replace
$|X_{n+1}-X_n|$ by $[X_{n+1}-X_n]^+$, i.e.\ to bound the
$p^{\rm th}$ moments of just the {\it positive part} of the increments.
Intuitively, this would correspond to allowing arbitrarily
large {\it negative} increments, and bounding only the large
{\it positive} increments.  The problem with this is that the
process is not sufficiently affected by its negative drift
when this is all concentrated into a few unlikely large jumps.
We give a counterexample to demonstrate this.

\itemm{}
Fix $0<\epsilon<1$, and consider the following {\it time-homogeneous} 
Markov chain $\{X_n\}$.  Let $X_0=0$, and for $n \ge 0$, let
$$
\P(X_{n+1} = 1 \, | \, X_n = 0) \ = \ 1
$$
$$
\P(X_{n+1} = x+1 \, | \, X_n = x > 0) \ = \ 1 - (1+\epsilon)/(x+1)
$$
$$
\P(X_{n+1} = 0 \, | \, X_n = x > 0) \ = \ (1+\epsilon)/(x+1)
$$

\itemm{}
Then (C1) is satisfied with $J=0$ and $a=\epsilon$.
Also, $\left[ X_{n+1}-X_n \right]^+ \le 1$, so (C2) would indeed
hold (for any $p>0$, and with $V=1$) if we replaced
$|X_{n+1}-X_n|$ by $[X_{n+1}-X_n]^+$.

\itemm{}
On the other hand, it is straightforward to see that
$\L(X_n)$ converges weakly to a stationary distribution $\pi(n)$,
which is such that $\pi(n) \sim C n^{-1-\epsilon}$ as $n\to\infty$.  
In particular, $\sum_n n \, \pi(n) =
\infty$.  It follows that $\E(X_n) \to \infty$, i.e.\ that $\{X_n\}$
is {\it not} stable in this case.  We conclude that Theorem~\ref{th:1}
does {\it not} continue to hold if we consider only the positive part of 
$X_{n+1}-X_n$ in condition (C2).  

\bigskip
\noindent \bf Remark. \rm
This last counter-example only works when $a \le V^{1/p}$.  In the
case where $X_n \geq 0$ for all $n$, this appears to be an
extremal counterexample, leading to the following open
question: 
\begin{quote}
Does Theorem~\ref{th:1} continues to hold 
for sufficiently large $a$ if we assume $X_n \geq 0$ and
replace $|X_{n+1}-X_n|$ by $[X_{n+1}-X_n]^+$ in (C2) ?
\end{quote}

Despite this counter-example, the hypotheses of Theorem~\ref{th:1} may
indeed be weakened to allow some large negative increments.  However, both
condition (C1) {\it and} condition (C2) must be identically modified
so that negative drift is still manifested.  This is the motivation for
having stated Corollary~\ref{th:9} as an extension to the main theorem.

\sect{Reduction to a martingale question}

We will derive Theorem~\ref{th:1} and Corollary~\ref{th:9} 
from the following martingale result.

\begin{th} \label{th:2}
Let $\{ M_n : n = 0 , 1 , 2 , \ldots \}$ be a sequence adapted to a 
filtration $\{ \F_n \}$ and let $\DD_n$ denote $M_{n+1} - M_n$.
Suppose that the sequence started at $M_1$ is a martingale (i.e., 
$\E (\DD_n \| F_n) = 0$ for $n \geq 1$), and that $M_0 \leq 0$.
Suppose further that for some $p > 2$ and $b > 0$ we have 
\begin{equation} \label{eq:moment}
\E (|\DD_n|^p \| \F_n) \leq b 
\end{equation}
for all $n$ including $n=0$.  Let $\tau = \inf \{ n > 0 : M_n \leq n \}$.  
Then for any $r \in (0 , p)$ there is a constant $C = C (b , p , r)$ 
such that 
\begin{equation} \label{eq:conclusion}
\E \left( (M_t^+)^r \one_{\tau > t} \right) \leq C t^{r - p} .
\end{equation}
\end{th}

We defer the proof of Theorem~\ref{th:2} until the following section.
In the remainder of this section, we assume Theorem~\ref{th:2}, and
derive Theorem~\ref{th:1} and Corollary~\ref{th:9} as consequences.

\begin{cor} \label{cor:1}
Let $\{ Y_n \}$ be adapted to $\{ \F_n \}$ with $Y_0 \leq 0$.
Suppose $\E (|\DE_n|^p \| \F_n) \leq B$ for all $n$ and
$\E (\DE_n \| \F_n) \leq 0$ for all $1 \leq n < \sigma$, where 
$\DE_n = Y_{n+1} - Y_n$ and $\sigma = \inf \{ n > 0 : Y_n \leq n \}$.  
Then for $0 < r < p$ there is a constant $K = K (B , p , r)$ such that
$$\E \left( (Y_t^+)^r \one_{\sigma > t} \right) \leq K t^{r - p} .$$
\end{cor}

\noindent{\sc Proof:} An easy fact useful here and later is that 
$z^+ \leq 1 + |z|^p$ and hence 
\begin{equation} \label{eq:p to 1}
\E |Z|^p \leq b \; \Rightarrow \; \E Z^+ \leq 1 + b .
\end{equation}
Recall (see e.g.\ Durrett 1996, p.~237) that the
supermartingale $\{ Y_{n \wedge \sigma} : n \geq 1 \}$ may be decomposed as
$Y_{n \wedge \sigma} = M_n - A_n$ where $\{ M_n : n \geq 1 \}$ is
a martingale and $\{ A_n : n \geq 1 \}$ is an increasing predictable
process with $A_1 = 0$.  Let $\mu_n$ denote $\E (\DE_n \| \F_n)$.
Then the increments $\DD_n := M_{n+1} - M_n$ satisfy 
$$\E (|\DD_n|^p \| \F_n) = \E (|\DE_n - \mu_n|^p \| \F_n) \leq 
   2^p \E (|\DE_n|^p + |\mu_n|^p \| \F_n) \leq 2^p (B + (1 + B)^p) .$$
Applying Theorem~\ref{th:2} to $\{ M_n \}$ with $b = 2^p (B + (1+B)^p)$
and $M_0 := Y_0$ yields
\begin{equation} \label{eq:C}
\E \left( (M_t^+)^r \one_{\tau > t} \right) \leq C t^{r-p} .
\end{equation}
When $\sigma > t$ it follows that $M_n \geq n + A_n$ for $1 \leq n \leq t$
and hence that $\tau > t$.  Also, when $\sigma > t$, we know that
$M_t = Y_t + A_t \geq Y_t$ and therefore that 
$$Y_t^+ \one_{\sigma > t} \leq M_t^+ \one_{\tau > t} .$$
The conclusion of the corollary now follows from~(\ref{eq:C}), with 
$K = C (2^p (B + (1+B)^p) , p , r)$.   $\Cox$

The above argument uses no properties of
the process $A_n$, other than its being nonincreasing and adapted.
In particular, it need not be predictable.  If the increments
$\DE_n$ can be decomposed into the sum of two parts, one 
satisfying the hypotheses of the corollary and one nonincreasing
and adapted, then the second of these can be absorbed into the
process $\{ A_n \}$ and the result will still hold.  Without
loss of generality, the second piece can be taken to be
$\DE_n \one_{\DE_n \leq Z_n}$ for some adapted nonpositive
$\{ Z_n \}$.  In other words, the moment condition need not
apply to the negative tail of the increment, as long as the
mean is still nonpositive when the negative tail is excluded.
We state this more precisely as the following corollary.

\begin{cor} \label{cor:better}
Let $\{ Y_n \}$ be adapted to $\{ \F_n \}$ with $Y_0 \leq 0$.
Let $\{ Z_n \}$ be any adapted nonpositive sequence.
Suppose $\E (|\DE_n|^p \one_{\DE_n > Z_n} \| \F_n) \leq B$ for all 
$n$ and $\E (\DE_n \one_{\DE_n > Z_n} \| \F_n) \leq 0$ for all 
$1 \leq n < \sigma$, where $\DE_n = Y_{n+1} - Y_n$ and 
$\sigma = \inf \{ n > 0 : Y_n \leq n \}$.  Then for $0 < r < p$,
$$\E \left( (Y_t^+)^r \one_{\sigma > t} \right) \leq K t^{r - p} .$$
$\Cox$
\end{cor}

We now use these corollaries to derive Theorem~\ref{th:1} and 
Corollary~\ref{th:9}.

\noindent{\sc Proof of Theorem~\protect{\ref{th:1}} from
Corollary~\protect{\ref{cor:1}}, and of Corollary~\ref{th:9} from
Corollary~\protect{\ref{cor:better}}}:  First assume that $a = 1$ and $J = 0$.
Given $\{ X_n \}$ as in the hypotheses of the theorem, fix an $N \geq 1$; 
we will compute an upper bound for $\E (X_N^+)^r$ that does not depend on $N$.  
Let $U := \max \{ k \leq N : X_k \leq 0 \}$ denote the last time up to 
$N$ that $X$ takes a nonpositive value.  Decompose according to the
value of $U$:
$$\E (X_N^+)^r = \sum_{k=0}^{N-1} \E \left( (X_N^+)^r \one_{U = k} \right) .$$

To evaluate the summand, define for any $k < N$ a process $\{ Y^{(k)}_n \}$
by $Y^{(k)}_n = (X_{k + n} + n) \one_{X_k \leq 0}$.  In other words, 
if $X_k > 0$ the process $\{ Y^{(k)}_N \}$ is constant at zero, and
otherwise it is the process $\{ X_n \}$ shifted by $k$ and compensated
by adding 1 each time step.  We apply Corollary~\ref{cor:1} to
the process $\{ Y^{(k)}_n \}$.  Hypothesis (C1) of Theorem~\ref{th:1},
together with the fact that $X_{k+j} > 0$ for $0 < j < \sigma^{(k)}$, 
imply that $\E (\DE_n \| \F_n) \leq 0$ when $1 \leq n \leq \sigma^{(k)}$.  
Also, $\E (|\DE_n|^p \| \F_n) \leq \E (|1 + X_{n+1} - X_n|^p \| \F_n) 
\leq B := 2^p (1 + V)$.  The conclusion is that 
$$\E \left( [(Y^{(k)}_{N-k})^+]^r \one_{\sigma^{(k)} > N - k} \right) \leq 
   K t^{r-p}$$  
with $K = K(V,p,r)$.
But for each $k$,
$$X_N^+ \one_{U = k} \leq Y^{(k)}_{N -  k} \one_{\sigma^{(k)} > N - k}$$
and it follows that 
$$\E \left( (X_N^+)^r \one_{U = k} \right) \leq K (N - k)^{r-p}.$$
Now sum to get
$$\E (X_N^+)^r \leq \sum_{k=0}^{N-1} K (N - k)^{r-p} \leq K \zeta (p-r) .$$
This completes the case $a = 1 , J = 0$.

For the general case, let $X_n' = (X_n - J) / a$.  This process is covered
by the analysis of the $a = 1 , J = 0$ case above, with $V / a^p$ 
in place of $V$.  We conclude that $\E (X_n')^r \leq c (p , 1 , V / a^2 , 
0 , r)$, and hence that $\E X_n \leq c (p , a , V , J , r) := 
J + a^r c(p , 1 , V / a^p , 0 , r)$.   

The proof of Corollary~\ref{th:9} from Corollary~\ref{cor:better} is
virtually identical.
$\Cox$

\sect{The proof of Theorem~\protect{\ref{th:2}}}

We now concern ourselves with the proof of Theorem~\ref{th:2}.  
We begin with two lemmas.

\begin{lem} \label{lem:Lp}
Let $\{ M_n \}$ be a martingale with $M_0 = 0$, and with increments bounded
in $L^p$:
$$\E (|M_n - M_{n-1}|^p \| \F_{n-1}) \leq L .$$
Then there is $c_p$ such that $\E |M_n|^p \leq c_p L n^{p/2}$.  
\end{lem}

\noindent{\sc Proof:}  Burkholder's inequality (see Stout 1974, Theorem~3.3.6;
Burkholder, 1966; Chow and Teicher, 1988, p.~396) tells us that 
for $p>1$, there is a constant $c_p$ for which
$$\E |M_n|^p \leq c_p \E \left ( \sum_{k=1}^n (M_k - M_{k-1})^2 
   \right)^{p/2}.$$
(In fact, for $p \ge 2$ we may take $c_p = (p-1)^p$, 
cf.\ Burkholder, 1988, Theorem 3.1.)
For any $Z_1 , \ldots , Z_n$, H\"older's inequality gives
$$\E |Z_1 + \cdots + Z_n|^{p/2} \leq n^{p/2} \max_{1 \leq k \leq n}
   \E |Z_k|^{p/2} .$$
Set $Z_k = (M_k - M_{k-1})^2$ and observe that $\E Z_k^{p/2} \leq 
\E ((\E (Z_k \| \F_{k-1})^{p/2}) \leq \E (\E |M_k - M_{k-1}|^p \| \F_{k-1}) 
\leq L$, so the conclusion of the lemma follows.    $\Cox$

\begin{lem} \label{lem:tau}
Assume the notation and hypotheses of Theorem~\ref{th:2}.
For $x > 0$, let $S_x = \inf \{k : M_k \geq x \}$ be the time to 
hit value $x$ or greater.  Then there is a $C' = C' (b , p)$ such that 
$$\P (\tau > S_x) \leq {C' \over x^{p/2}} \, .$$ 
\end{lem}

\noindent{\sc Proof:}  Fix $x \ge 1$ and bound in two ways the quantity 
$\E \left | M_{\tau \wedge S_x} \right |^p$.  
First, since $\{ M_{\tau \wedge S_x \wedge n} : n \geq 1 \}$ is a 
martingale, $|x|^p$ is convex, and $\tau \wedge S_x \geq 1$ is a 
stopping time bounded above by $x$, we have 
\begin{equation} \label{eq:upside}
\E \left | M_{\tau \wedge S_x} \right |^p \leq \E \left | M_x \right |^p .
\end{equation}
Using Lemma~\ref{lem:Lp} gives $\E |M_x - M_1|^p \leq c_p b x^{p/2}$,
and since $\E |M_1|^p \leq b$, this yields
\begin{equation} \label{eq:upside2}
\E |M_x|^p = ||M_x||_p^p \leq (||M_1||_p + ||M_x - M_1||_p)^p 
   \leq \left ( b^{1/p} + (c_p b x^{p/2})^{1/p} \right )^p 
   \leq c'(p,b) x^{p/2} 
\end{equation}
with $c'(p,b) := c_p b (1 + c_p^{-1})^p$.  On the other hand, 
on the event $\{\tau > S_x\}$ we have $M_{\tau \wedge S_x}
= M_{S_x} \ge x$, so that
$$x^p \P (\tau > S_x) \leq \E |M_{\tau \wedge S_x}|^p , $$ 
and combining this with~(\ref{eq:upside}) and~(\ref{eq:upside2}) gives
$$\P (\tau > S_x) \leq x^{-p} c'(p,b) x^{p/2} $$
which proves the result for $x \ge 1$.  Finally, for $x < 1$ we use
$\P( \tau > S_x ) \le 1$, so the lemma follows with
$C'(p,b) := \max(1, \, c'(p,b))$.
$\Cox$ 

\noindent{\sc Proof of Theorem}~\ref{th:2}: Let $T = \inf \{ k \geq 0 : 
\DD_k \geq t/4 \}$ be the time of the first large jump.  Since
$\tau > t$ implies $S_x < x$ for all $x \leq t$, we can write
\begin{equation} \label{eq:star}
\E \left( (M_t^+)^r \one_{\tau > t} \right) 
= \E \left( (M_t^+)^r \one_G \right) 
+ \E \left( (M_t^+)^r \one_H \right) ,
\end{equation}
where $G = \{ T \geq S_{t/2} < t < \tau \}$ and $H = \{ T < S_{t/2} 
< t < \tau \}$.  

To bound the first term, abbreviate $S := S_{t/2}$ and begin by 
observing that $M_S \leq 3t/4$ on $G$, since the level $t/2$ or higher 
has just been obtained and the increment was no more than $t/4$.  Thus 
$$\E \left( (M_t^+)^r \one_G \right) 
\leq \P (\tau > S) \E ((M_t^+)^r \one_{\tau > t} \| \F_S) . $$
The first factor may be bounded via Lemma~\ref{lem:tau}:
\begin{equation} \label{eq:S}
\P (\tau > S) \leq {2^{p/2} C' \over t^{p/2}} \, .
\end{equation}
The second factor is bounded using the formula
\begin{equation} \label{eq:parts}
\E \left( Z^r \one_{Z > u} \right) = u^r \P (Z > u) + \int_u^\infty r y^{r-1} 
   \P (Z > y) \, dy .
\end{equation}
By Lemma~\ref{lem:Lp} conditionally on $\F_S$, $\E (|M_t - M_S|^p \| F_S) 
\leq c_p b (t-S)^{p/2} \leq c_p b t^{p/2}$.  Hence by Markov's inequality,
$\P (M_t - M_S \geq y \| \F_S) \leq c_p b t^{p/2} / y^p$.  Therefore,
\begin{eqnarray*}
\E ((M_t^+)^r \one_{\tau > t} \| \F_S) & \leq & \E ((M_t^+)^r \one_{M_t > t} 
   \| \F_S) \\
& = & t^r \P (M_t \geq t \| \F_S) + \int_t^\infty r y^{r-1} 
   \P (M_t \geq y \| \F_S) \, dy \\
& \leq & t^r \P (M_t - M_S \geq t/4 \| \F_S) + \int_{t/4}^\infty 
   r y^{r-1} \P (M_t - M_S \geq y \| \F_S) \, dy \\
& \leq & c_p b 4^p t^{r - p/2} + \int_{t/4}^\infty 
   r y^{r-1} c_b p t^{p/2} y^{-p} \, dy \\
& \leq & c_p b (4^p + {r \over p - r} 4^{p-r}) t^{r-p/2} \\
& \leq & c_2 (b , p , r) t^{r - p/2} ,
\end{eqnarray*}
where $c_2 (b , p , r) := c_p b (4^p + 4^{p-r} r / (p-r))$.  Combining
with~(\ref{eq:S}) gives
\begin{equation} \label{eq:term1}
\E \left( (M_t^+)^r \one_G \right) \leq 2^{p/2} c_2 C' t^{r - p} .
\end{equation}

We will bound the second term by decomposing according to the value 
of $T$.  A preliminary computation is to bound the quantity
$\E ((M_t^+)^r \one_{T = k , M_t > t} \| \F_k)$.  Break this into three
pieces: the part up to time $k$, the jump at time $k$, and the
part from time $k+1$ to time $t$.  For any $0 < r < p-1$, 
$|x+y+z|^r \leq 3^r (|x|^r + |y|^r + |z|^r)$ (use convexity
when $r \geq 1$ and sublinearity when $r \leq 1$).  Hence
$$(M_t^+)^r \leq 3^r \left[ (M_k^+)^r + (\DD_k^+)^r + ((M_t - M_{k+1})^+)^r 
\right] .$$
The event $\{ T = k \}$ implies $M_k \leq 3t/4$, and is also in the 
initital $\sigma$-field of the martingale $\{ M_n - M_{k+1} : 
n \geq k+1 \}$.  Therefore, when we condition on $\F_k$, we get
\begin{eqnarray*}
\E ((M_t^+)^r \one_{T = k , M_t > t} \| \F_k) & \leq & 3^r \left [ 
   \left ( {3t \over 4} \right )^r \P (\DD_k \geq t/4 \| \F_k) \right. \\
&& + \E ((\DD_k^+)^r \one_{\DD_k \geq t/4} \| \F_k) \\
&& + \left. \E (|M_t - M_{k+1}|^r \one_{T = k , M_t - M_{k+1} \geq t/4} 
   \| \F_k) \right ] 
\end{eqnarray*}

The moment condition $\E (|\DD_k|^p \| \F_k) \leq b$ implies that
$\P (\DD_k \geq y) \leq b y^{-p}$, hence the first of these
contributions is at most 
$$3^r ({3t \over 4})^r b ({t \over 4})^{-p} .$$
Using~(\ref{eq:parts}) again, we bound the second of the three
contributions by
$$3^r ({t \over 4})^r \P (\DD_k \geq {t \over 4} \| \F_k) 
   + 3^r \int_{t/4}^\infty r y^{r-1} \P (\DD_k \geq y) \, dy$$
which is at most 
$$3^r b ({t \over 4})^{r-p} + 3^r {r \over p-r} ({t \over 4})^{r-p} .$$
Lemma~\ref{lem:Lp} implies $\E (|M_t - M_{k+1}|^r \| \F_{k+1}) \leq 
c_p b t^{r/2}$, while $\one_{T=k} \in \F_{k+1}$ and has conditional
expectation at most $b (t/4)^{-p}$ given $\F_k$.  Therefore the
third contribution is bounded by 
$$3^r b ({t \over 4})^{-p} c_p b t^{r/2} .$$
Summing the three contributions gives
\begin{equation} \label{eq:ub3}
\E ((M_t^+)^r \one_{T = k , M_t > t} \| \F_k) \leq c_3 t^{r-p}
\end{equation}
where $c_3 = 3^r 4^p b (c_p b + {p \over p-r} + 3^r)$. 

Now we bound the second term, by 
decomposing according to the value of $T$.
\begin{eqnarray}
\E \left( (M_t^+)^r \one_H \right) &
= & \sum_{k=0}^{\lfloor t/2 \rfloor} 
\E \left( (M_t^+)^r \one_H \one_{T = k} \right)
   \nonumber \\
& = & \sum_{k=0}^{\lfloor t/2 \rfloor} \E \left [ \E ((M_t^+)^r \one_H \one_{T = k} \| \F_k)
   \right ] \label{eq:all k} .
\end{eqnarray}
The event $\{ \tau > k \}$ is in $\F_k$ and contains the event
$H \cap \{ T = k \}$, so we have
$$ \E ((M_t^+)^r \one_H \one_{T = k} \| F_k) \leq \one_{\tau > k} 
   \, \E ((M_t^+)^r \one_{T = k} \| \F_k)$$ 
and hence
$$ \E ((M_t^+)^r \one_H) \leq \sum_{k=0}^{\lfloor t/2 \rfloor} \P (\tau > k) 
   \, \E \left [ \E ((M_t^+)^r) \one_{T = k} \| \F_k) \right ] . $$
Plugging in the upper bound~(\ref{eq:ub3}) and using Lemma~\ref{lem:tau} 
gives
$$\E ((M_t^+)^r \one_H) \leq \sum_{k=0}^{\lfloor t/2 \rfloor} C' k^{-p/2} c_3 t^{r-p}$$
and summing yields a bound of
\begin{equation} \label{eq:term2}
\E ((M_t^+)^r \one_H) \leq c_4 t^{r-p}
\end{equation}
for the second term,
where $c_4 := C' c_3 \zeta(p/2)$.  
By~(\ref{eq:star}), the two bounds~(\ref{eq:term1}) and~(\ref{eq:term2}) 
together imply the conclusion of Theorem~\ref{th:2}.   $\Cox$

This completes the proof of Theorem~\ref{th:2}, and hence also the
proof of Theorem~\ref{th:1} and Corollary~\ref{th:9}.

\bigskip
\noindent \bf Acknowledgements. \rm
We are grateful to Allan Borodin for bringing this problem to our
attention.  We thank Jim Fill, Leslie Goldberg, Jim Kuelbs, Tom Kurtz,
Jeremy Quastel, Tom Salisbury, and Richard Tweedie for very helpful comments.

\end{document}